\documentclass[a4paper]{amsart}
\usepackage{amsmath,amsfonts,amssymb,amsthm}
\usepackage{comment}
\usepackage{hyperref}

\newtheorem{theorem}{Theorem}
\newtheorem{corollary}[theorem]{Corollary}
\newtheorem{lemma}[theorem]{Lemma}
\newtheorem{proposition}[theorem]{Proposition}
\newtheorem{definition}{Definition}
\theoremstyle{remark}
\newtheorem{remark}{Remark}

\def\Z{{\mathbb Z}}
\def\Q{{\mathbb Q}}

\def\C{{\mathbb C}}
\def\P{{\mathbb P}}
\def\F{{\mathbb F}}
\def\p{{\mathfrak p}}

\DeclareMathOperator{\ch}{char}

\DeclareMathOperator{\disc}{disc}
\DeclareMathOperator{\ord}{ord}
\DeclareMathOperator{\diag}{diag}

\DeclareMathOperator{\Br}{Br}
\DeclareMathOperator{\GL}{GL}
\DeclareMathOperator{\SL}{SL}
\DeclareMathOperator{\PGL}{PGL}
\DeclareMathOperator{\N}{N}
\DeclareMathOperator{\T}{Tr}
\def\NLK{\N_{L/K}}
\def\TLK{\T_{L/K}}
\def\Ks{K^*}
\def\Ksq{K^*{}^2}
\def\Ksc{K^*{}^3}
\def\Ksmc{\Ks/\Ksc}
\def\Ls{L^*}
\def\Lsc{L^*{}^3}
\def\Lsmc{\Ls/\Lsc}

\def\D{\Delta}

\newcommand{\BC}{{\mathcal {BC}}}
\renewcommand{\AA}{{\mathcal{A}}}
\newcommand{\CC}{{\mathcal {C}}}

\begin{document}
\title{On the equivalence of binary cubic forms}
\author{J. E. Cremona}
\address{Mathematics Institute, University of Warwick, Coventry, CV4
  7AL, UK.}
\email{J.E.Cremona@warwick.ac.uk}
\date{\today}

\begin {abstract}
We consider the question of determining whether two binary cubic forms
over an arbitrary field~$K$ whose characteristic is not~$2$ or~$3$ are
equivalent under the actions of either~$\GL(2,K)$ or~$\SL(2,K)$,
deriving two necessary and sufficient criteria for such equivalence in
each case. One of these involves an algebraic invariant of binary
cubic forms, which we call the Cardano invariant, as it is closely
connected to classical formulas; it also appears in the work of
Bhargava {\it et al.}. The second criterion is in terms of
the base field itself, and also gives explicit matrices in~$\SL(2,K)$
or~$\GL(2,K)$ transforming one cubic into the other, if any exist, in
terms of the coefficients of bilinear factors of a bicovariant of the
two cubics.  We also consider automorphisms of a single binary cubic
form, show how to use our results to test equivalence of binary cubic
forms over an integral domain such as~$\Z$, and briefly recall some
connections between binary cubic forms and the arithmetic of elliptic
curves.

The methods used are elementary, and similar to those used in our
work with Fisher concerning equivalences between binary
quartic forms.
\end {abstract}

\maketitle

\section{Introduction}
\label{sec:intro}
We consider binary cubic forms over an arbitrary field~$K$ whose
characteristic is not~$2$ or~$3$, and establish two necessary and
sufficient conditions under which two cubics are equivalent under the
actions of~$\GL(2,K)$ or~$\SL(2,K)$, which yield tests which are
simple to apply. One of these conditions involves an algebraic
invariant of binary cubic forms, which we call the Cardano invariant,
as it is closely connected to classical formulas.  A more general
version of this invariant, involving ideal classes in quadratic rings,
appears in the work of Bhargava, Elkies and Shnidman
\cite{bhargava-elkies-shnidman2020}, in the classification of
$\SL(2,D)$-orbits of binary cubic forms over a Dedekind domain~$D$:
see Theorem~18 and Corollary~20
of~\cite{bhargava-elkies-shnidman2020}.  The invariant was first
introduced by Bhargava, in the case where the base ring is~$\Z$,
in~\cite{bhargava2004}; it is defined as an element of the quadratic
resolvent algebra associated to the cubic form.  Our second result
gives a test for equivalence in terms of the base field
itself, which also gives explicit matrices transforming one cubic into a
second when they are equivalent.

We also consider automorphisms of a single binary cubic form,
recovering a result of Xiao in~\cite{xiao2019binary}, which also
follows (for $K=\Q$) from the Delone-Faddeev parametrization, and we
show how to use our results to test equivalence of binary cubic forms
over an integral domain such as~$\Z$.  Finally, we make some remarks
on how such results relate to the arithmetic of elliptic curves, as
in~\cite{bhargava-elkies-shnidman2020}.

The methods used here, which are all elementary, are similar to those
used in our work \cite{CremonaFisherQuartics} with Fisher, concerning
equivalences between binary quartic forms.  While our first criterion
for~$\SL(2,K)$-equivalence may be found in the literature, our account
is self-contained and more elementary, compares~$\SL(2,K)$-equivalence
with ~$\GL(2,K)$-equivalence, and our second criterion has the merit
of being purely algebraic, without a need to extend the base field,
together with the benefit of giving explicit transformation matrices.

Over an algebraically closed field~$\overline{K}$, all binary cubic
forms with nonzero discriminant are $\GL(2,\overline{K})$-equivalent
(to $XY(X+Y)$, for example); the issue we address is whether
equivalences exist which are defined over the ground field~$K$ itself,
and if so, to find them explicitly using only arithmetic in~$K$,
without having to extend to the splitting field.  We also consider the
question of~$\SL(2,K)$-equivalence.

In order to state our results, we make some definitions.  Let~$\BC(K)$
denote the set of all binary cubic forms in~$K[X,Y]$ with non-zero
discriminant, and~$\BC(K;\D)$ to be the subset with discriminant~$\D$,
for each~$\D\in\Ks$.  We define the \emph{resolvent algebra}~$L$
associated to~$\D\in\Ks$ to be the quadratic \'etale algebra
\[
L = K[\delta] = K[X]/(X^2+3\D).
\]
We define the \emph{twisted action} of~$\GL(2,K)$ on~$\BC(K)$ as
follows: for $g\in\BC(K)$ and
$M=\begin{pmatrix}r&s\\t&u\end{pmatrix}\in\GL(2,K)$, set
\[
g^M(X,Y) = \det(M)^{-1} g(r X + t Y, s X + u Y) = \det(M)^{-1} g(X',Y'),
\]
where $(X'\ Y')=(X\ Y)M$.  Both the untwisted action (without the
determinant factor) and this action appear in the literature: for
example, the untwisted action is used
in~\cite{bhargava-elkies-shnidman2020} and is also the subject
of~\cite{kulkarni-ure2021}, while the twisted action is used
in~\cite{bhargava-shankar-tsimerman2013}.  Clearly they are the same
for~$M\in\SL(2,K)$, and the difference is only minor for our purposes.

In Section~\ref{sec:cardano} we define the \emph{Cardano
invariant}~$z(g)$ for~$g\in\BC(K;\D)$ to be an element of~$\Lsmc$, and
show that it has the following properties, where, as
in~\cite{bhargava-elkies-shnidman2020}, we denote the kernel of the
induced norm map $\NLK: \Lsmc \to \Ksmc$ by $(\Lsmc)_{N=1}$.  We name
this group the \emph{Cardano group} for~$\D$.
\begin{theorem}\label{thm:main1}
  Let $K$ be any field with~$\ch(K)\not=2,3$.
  Let~$\D\in\Ks$, let~$L$ be the resolvent
  algebra~$K[X]/(X^2+3\D)$, and let~$z\colon \BC(K;\D)\to\Lsmc$ be
  the Cardano invariant map.
  \begin{enumerate}
  \item $z(g) \in (\Lsmc)_{N=1}$ for all~$g\in\BC(K;\D)$;
  \item $z(g)=1$ if and only if~$g$ is reducible over $K$;
  \item $g_1,g_2\in \BC(K;\D)$ are $\SL(2,K)$-equivalent if and only
    if $z(g_1)=z(g_2)$;
  \item $g_1,g_2\in \BC(K;\D)$ are $\GL(2,K)$-equivalent if and only
    if $z(g_1)=z(g_2)^{\pm1}$ (equivalently, if and only if $z(g_1)$ and~$z(g_2)$
    generate the same subgroup of the Cardano group $(\Lsmc)_{N=1}$);
  \item $z$ induces bijections between the $\SL(2,K)$-orbits on
    $\BC(K;\D)$ and the Cardano group, and between the
    $\GL(2,K)$-orbits on~$\BC(K;\D)$ and its cyclic subgroups.
  \end{enumerate}
\end{theorem}
In fact, discriminant-preserving transformations all have
determinant~$\pm1$ (see Proposition~\ref{prop:SL2GL2} below); the
Cardano invariant is preserved by those with determinant~$+1$, and
inverted by those of determinant~$-1$.

The analogue of part~(3) of this theorem for binary cubic forms
over~$\Z$ follows from Theorem~13 of~\cite{bhargava2004}, and the
proof there carries over without significant change: our Cardano
invariant is denoted~$\delta$ in~\cite{bhargava2004}; a quadratic
ring~$S$ (over~$\Z$) is fixed whereas we fix the
discriminant~$\D\in\Ks$, and the ideal $I$ of \cite{bhargava2004} is
irrelevant, as our base ring is a field.  This assumption
enables us to give a simpler proof than in~\cite{bhargava2004}.

The last part of the Theorem, for the case of~$\SL(2,K)$-orbits, is
the same as Corollary~20 of~\cite{bhargava-elkies-shnidman2020}.

This theorem implies that the $\SL(2,K)$-orbits on~$\BC(K;\D)$ carry a
group structure, in which the identity is the class of reducible
cubics.  Testing whether~$z(g_1)=z(g_2)$ amounts to testing whether a
certain monic cubic over~$K$, constructed from~$g_1$ and~$g_2$, has a
root in~$K$. In terms of the group structure, we are testing that
the~$\SL(2,K)$-orbits of~$[g_1]$ and~$[g_2]$ are equal by testing
whether~$[g_1][g_2]^{-1}$ is trivial.  For an interpretation of this
group in terms of Galois cohomology, see
Section~\ref{sec:elliptic-curves}.

The reason we call $z(g)$ the Cardano invariant is that the classical
formula of Cardano for solving cubic equations involves the cube root
of an expression in $K(\sqrt{-3\D})$, and for monic cubics this
expression is precisely~$z(g)$.  See Remark~\ref{rmk:Cardano} below
for a precise statement.

Our second result gives conditions for two cubics with the same
discriminant to be~$\SL(2,K)$-{} or~$\GL(2,K)$-equivalent, which only
involve factorization of multivariate polynomials over~$K$,
and which also provide explicit matrices~$M\in\GL(2,K)$ such that
$g_1^M=g_2$, when such matrices exist.  For brevity, we only state
here the result for~$\SL(2,K)$; see Proposition~\ref{prop:11factors}
and Theorem~\ref{thm:equiv2} for the general case.

For two cubic forms $g_1,g_2$ over~$K$, with cubic covariants
$G_1,G_2$ respectively (whose definition we recall in
Section~\ref{sec:basics}), we define the \emph{bicubic bicovariant}
$B_{g_1,g_2} \in K[X_1,Y_1,X_2,Y_2]$ to be
\[
B_{g_1,g_2}(X_1,Y_1,X_2,Y_2) = g_2(X_2,Y_2)G_1(X_1,Y_1) - G_2(X_2,Y_2)g_1(X_1,Y_1);
\]
this is a bihomogeneous form of degree~$3$ in each pair of variables
$X_1,Y_1$ and $X_2,Y_2$.

\begin{theorem}\label{thm:main2}
  Let $g_1,g_2$ be binary cubic forms over~$K$ with
  $\disc(g_1)=\disc(g_2)\not=0$ and with bicovariant $B_{g_1,g_2}$.
  Then the following are equivalent:
  \begin{enumerate}
  \item
    $B_{g_1,g_2}$ has a bilinear factor in $K[X_1,Y_1,X_2,Y_2]$;
  \item
    $g_1$ and~$g_2$ are $\SL(2,K)$-equivalent.
  \end{enumerate}
  More precisely, $B_{g_1,g_2}$ has the factor
  $-s X_1X_2 + r X_1Y_2 - u Y_1X_2 + t Y_1Y_2$
  (up to scaling) if and only if $g_1=g_2^M$, where
  $M=\begin{pmatrix}r&s\\t&u \end{pmatrix}\in\SL(2,K)$.
\end{theorem}

\begin{remark} Over the algebraic closure~$\overline{K}$, the
bicovariant~$B_{g_1,g_2}$ factors as a product of three bilinear
factors (that is, bihomogeneous factors linear in each pair of
variables), and there are precisely three
matrices~$M\in\SL(2,\overline{K})$ with $g_1=g_2^M$. If~$M$ is one of
these, then the other two are $MT$, $MT^2$
where~$T\in\SL(2,\overline{K})$ satisfies $T^3=I$ and $g_1^T=g_1$.
Over~$K$ itself, there can be at most one~$M\in\SL(2,K)$
satisfying~$g_1^M=g_2$, unless $\D$ is a square in~$K$, when there are
three such (if any), related as above.
\end{remark}

\begin{remark}
We show in Section~\ref{sec:equivalence} below (see
Proposition~\ref{prop:SL2GL2}) that $g_1,g_2\in\BC(K;\D)$ are
$\GL(2,K)$-equivalent if and only if $g_1(X,Y)$ is
$\SL(2,K)$-equivalent to~$g_2(X,Y)$ or to~$g_2(X,-Y)$, the matrix
transforming~$g_1$ into~$g_2$ having determinant~$+1$ or~$-1$
respectively.  Transformations with determinant~$-1$ are associated to
bilinear factors of a twisted form of~$B_{g_1,g_2}$; see see
Proposition~\ref{prop:11factors}.
\end{remark}

In the following section we recall the definitions of the invariants,
seminvariants and covariants associated to a binary cubic form. In
Section~\ref{sec:cardano} we define the Cardano covariant and Cardano
invariant and establish their basic properties.  In
Section~\ref{sec:equivalence} we prove both of our main results. In
the last three sections,
we discuss automorphisms of binary cubic forms, how to extend our
results to integral forms, and the connections between binary cubic
forms and the arithmetic of elliptic curves.

\subsubsection*{Acknowledgements} We thank Tom Fisher and a referee for
helpful comments and some additional references.

\section{Binary cubics, their invariants and covariants}
\label{sec:basics}

Let $\BC(K)$ denote the space of binary cubic forms in~$K[X,Y]$ with
nonzero discriminant, and for each~$\D\in\Ks$, let $\BC(K;\D)$ be the
subset of those forms with discriminant~$\D$.  Recall that for $g(X,Y)
= aX^3+bX^2Y+cXY^2+dY^3\in\BC(K)$, the \emph{discriminant}~$\disc(g)$
of~$g$ is given by
\[
   \D = b^2c^2-4ac^3-4b^3d-27a^2d^2+18abcd;
\]
we also define the seminvariants\footnote{Seminvariants are the
leading coefficients of covariants; they are invariant under upper
triangular transformations, and include invariants as a special case.}
\[
P = b^2-3ac,   \qquad\hbox{and}\qquad
U = 2b^3+27a^2d-9abc,
\]
which satisfy the syzygy
\[
4P^3 = U^2 + 27\D a^2.
\]
These are the leading coefficients of two covariant binary forms, the
quadratic \emph{Hessian}
\[
H(X,Y) = (b^2-3ac)X^2 + (bc-9ad)XY + (c^2-3bd)Y^2,
\]
which has discriminant~$-3\D$, and the cubic
\begin{multline*}
G(X,Y)  = (2b^3+27a^2d-9abc)X^3 + 3(b^2c+9abd-6ac^2)X^2Y \\ -3(bc^2+9acd-6b^2d)XY^2 - (2c^3+27ad^2-9bcd)Y^3,
\end{multline*}
which has discriminant~$729\D^3$; these satisfy the syzygy
\begin{equation}\label{eqn:syzygy}
4H(X,Y)^3 = G(X,Y)^2 + 27\D g(X,Y)^2,
\end{equation}
extending the seminvariant syzygy which is recovered on
setting~$(X,Y)=(1,0)$.

We define the \emph{twisted action} of~$\GL(2,K)$ on~$\BC(K)$ as
follows: for $g\in\BC(K)$ and $M=\begin{pmatrix}r&s\\t&u\end{pmatrix}\in\GL(2,K)$, set
\[
g^M(X,Y) = \det(M)^{-1} g(r X + t Y, s X + u Y) = \det(M)^{-1} g(X',Y'),
\]
where $(X'\ Y')=(X\ Y)M$.  Both the untwisted action (without the
determinant factor) and this action appear in the literature: for
example, the untwisted action is used
in~\cite{bhargava-elkies-shnidman2020} and is also the subject
of~\cite{kulkarni-ure2021}, while the twisted action is used
in~\cite{bhargava-shankar-tsimerman2013}.  Clearly they are the same
for~$M\in\SL(2,K)$, and the difference is only minor for our purposes
since we will mostly be concerned with matrices with
determinant~$\pm1$, and $-g = g^M$ with $M=-I$.

The discriminant and other covariants of~$g$ and~$g^M$ are related
by\footnote{The action on the binary quadratic form $H(g)$ is untwisted:
$H(g)^M(X,Y)=H(g)(X',Y')$.  The untwisted action on cubics gives
a factor of $\det(M)^k$ on the right with $k=6,2,3$ respectively.}
\begin{align*}
\D(g^M) &= \det(M)^2\D(g);\\
H(g^M)  &= H(g)^M;\\
G(g^M)  &= \det(M) G(g)^M;
\end{align*}
indeed, the definition of covariants is precisely that such identities
hold; the invariant~$\D$ is just a covariant of degree~$0$.

Below we will also need that $g\in\BC(K)$ is coprime to its
covariant~$G$, their resultant being~$8\D^3$.

\section{The resolvent algebra and the Cardano invariant}
\label{sec:cardano}

To each~$\D\in\Ks$, we associate the \emph{resolvent algebra}
\[
L = K[\delta] = K[T]/(T^2+3\D),
\]
a quadratic \'etale algebra over~$K$, with a fixed generator~$\delta$
(the image of~$T$ in the quotient) such that $\delta^2=-3\D$; it is a
field unless $-3\D$ is a square in~$K$.  We denote by a bar the
nontrivial $K$-automorphism of~$L$, mapping~$u+v\delta\mapsto
u-v\delta$.  Let~$L^*$ denote the unit group of~$L$, which consists of
the elements with nonzero norm; the norm map $\NLK: L \to K$ maps
$z=u+v\delta \mapsto z\overline{z}=u^2+3\D v^2$. We extend these to
maps $L[X,Y]\to L[X,Y]$ and $L[X,Y]\to K[X,Y]$, and also use the same
notation for the induced group homomorphisms~$\Lsmc\to\Lsmc$
and~$\Lsmc\to\Ksmc$.

For~$M\in\GL(2,K)$ and~$g\in\BC(K)$, since~$\D(g^M)=\det(M)^2\D(g)$,
the resolvent algebras of~$g$ and~$g^M$ are isomorphic via
$\delta\mapsto\det(M)\delta$.  Transforms which preserve the
discriminant therefore have determinant~$\pm1$; however,
we regard~$\delta$ as a fixed generator
of~$L$, associated to~$\Delta$, which is the same
for all cubics in~$\BC(K;\D)$.

\begin{definition}
The \emph{Cardano covariant} $C(X,Y)$ of $g\in\BC(K;\D)$ is
\[
C(X,Y) = \frac{1}{2}(G(X,Y)+3\delta g(X,Y)) \in L[X,Y].
\]
\end{definition}
Classically, this is called an ``irrational'' or algebraic covariant
because its coefficients lie in the extension~$L$ rather than in $K$
itself.  It is a covariant for the action of~$\SL(2,K)$; in general,
we have
\begin{equation}\label{eq:Cardano-transform}
C(g^M) = \frac{1}{2}(\det(M)G + 3\delta g)^M;
\end{equation}
this would be equal to~$\det(M)C(g)^M$ if the value of~$\delta$ for
discriminant~$\det(M)^2\D$ were fixed to be $\det(M)\delta$, but it is
not possible to do this consistently for transformations with
both determinants~$\pm\det(M)$.

At first, we will restrict our attention to the action of~$\SL(2,K)$.

In terms of~$C$, the covariant syzygy may be simply written
\[
\NLK(C) = C\overline{C} = H^3.
\]
Hence, for all $x,y\in K$, not both~$0$, except for those
satisfying~$H(x,y)=0$ (which only exist when $-3\D$ is a square
in~$K$), the value $C(x,y)$ lies in~$L^*$ and satisfies
$\NLK(C(x,y))=H(x,y)^3$ with $H(x,y)\in K^*$.

We will define the Cardano invariant in terms of values of the Cardano
covariant.  Before we give the general definition, first consider
cubics with~$P=H(1,0)\not=0$.  Set $z = C(1,0) =
\frac{1}{2}(U+3a\delta)$; then $\NLK(z)=P^3 \in \Ksc$, so $z\in L^*$,
and a provisional definition of the Cardano invariant of~$g$ with
$P\not=0$ is simply this quantity~$z$.  One can show (see
Proposition~\ref{prop:linear_factor} below) that~$z\in\Lsc$ if and
only if~$g$ has a linear factor in~$K[X,Y]$, and that for two
cubics~$g_1,g_2$ with the same discriminant~$\D$, and hence the same
resolvent algebra, they are $\SL(2,K)$-equivalent if and only if their
$z$-invariants are equal modulo~$\Lsc$.  Instead, however, we proceed
as follows, leading to a general definition of the Cardano
invariant~$z(g)$ as a well-defined element of~$\Lsmc$ for
all~$g\in\BC(K;\D)$.

Observe that~$C$, whose norm~$H^3$ is a cube in~$K[X,Y]$, is itself
the cube of a linear form in~$L[X,Y]$, up to a constant factor.
Explicitly, we have
\[
27 c_0^2 C(X,Y) = (3 c_0 X + c_1 Y)^3,
\]
where $C(X,Y) = c_0X^3+c_1X^2Y+c_2XY^2+c_3Y^3$, so $c_0=z$ (as defined
above) and~$c_1=\frac{3}{2}(b^2 c - 6 a c^2 + 9 a b d + b \delta)$.
(This algebraic identity may be readily checked directly, or derived
by writing $C$ as a constant times the cube of a linear form, and
differentiating twice.) This already shows that if $P\not=0$, so
that~$c_0\in \Ls$, then for every $x,y\in K$ such that $H(x,y)\not=0$,
we have $C(x,y)\in \Ls$ and~$C(x,y)\equiv c_0\mod{\Lsc}$.  For the
general case, we use the following identity:

\begin{proposition}\label{prop:Cmodcubes}
The identity
\[
   C(X_1,Y_1)^2 C(X_2,Y_2) = F(X_1,Y_1,X_2,Y_2)^3
\]
holds, where $F\in L[X_1,Y_1,X_2,Y_2]$ is given by
\[
3F(X_1,Y_1,X_2,Y_2) =
(3c_0X_1^2 + 2c_1X_1Y_1 + c_2Y_1^2) X_2 +
(c_1X_1^2 + 2c_2X_1Y_1 + 3c_3Y_1^2) Y_2.
\]
\end{proposition}

\begin{proof}
This identity reduces to the previous one on specialising
$(X_1,Y_1)=(1,0)$, and may be checked using computer algebra.
\end{proof}

Also, $\NLK F(X_1,Y_1,X_2,Y_2) = H(X_1,Y_1)^2 H(X_2,Y_2)$ and~$C(X,Y)
= F(X,Y,X,Y)$.

\begin{corollary}
The value of $C(x,y)\in \Lsmc$ is independent of $(x,y)\in K\times K$,
provided that $H(x,y)\not=0$, so that~$C(x,y)$ is a unit.
\end{corollary}
\begin{proof}
By the proposition, for $x_1,y_1,x_2,y_2$ such that
$H(x_1,y_1),H(x_2,y_2)\not=0$, we have $C(x_1,y_1), C(x_2,y_2) \in
\Ls$ and, modulo~$\Lsc$, we have
\[
C(x_1,y_1)^{-1}C(x_2,y_2) \equiv C(x_1,y_1)^2C(x_2,y_2) \equiv
F(x_1,y_1,x_2,y_2)^3 \equiv1.
\]
\end{proof}
Hence we may define the {\em Cardano invariant} $z(g)$ of $g\in\BC(K;\D)$ as
an element of~$\Lsmc$, by setting
\[
z(g) = C(x,y) \in \Lsmc,
\]
for any choice of~$x,y\in K$ such that~$C(x,y)$ is a unit.  This is
well-defined by the corollary, and if $P=H(1,0)\not=0$ then we may
take $z(g)=C(1,0)=z$ (modulo cubes), as above.  In all cases we have,
modulo~$\Ksc$,
\[
\NLK(z(g)) \equiv \NLK(C(x,y)) \equiv H(x,y)^3 \equiv 1,
\]
so that $z(g) \in (\Lsmc)_{N=1} := \ker(\NLK: \Lsmc \to \Ksmc)$, the
\emph{Cardano group}.  This establishes the first part of
Theorem~\ref{thm:main1}.  The second part is the following:

\begin{proposition} \label{prop:linear_factor}
$g\in\BC(K)$ has a linear factor in $K[X,Y]$ if and only if $z(g)=1$.
\end{proposition}
\begin{proof}
Writing $g=aX^3+bX^2Y+cXY^2+dY^3$, we first observe that when $a=0$ we
have $P=b^2\not=0$ (as otherwise $\Delta=0$), so
$z(g)=z=U/2=b^3\in\Ksc$, and $g$ has the factor~$Y$.

Now assume that $a\not=0$.  The equations
$z=\frac{1}{2}(U+3a\delta)=(x+y\delta)^3$, together with\footnote{If
$K$ contains cube roots of unity, we may scale $x$ and $y$ by a cube
root of unity if necessary, so that $\NLK(x+y\delta)=P$ exactly.}
$\NLK(x+y\delta)=P$, have a solution~$x,y\in K$ if and only if
\[
f(x) = 0 \quad\text{and}\quad y = 9a / f'(x),
\]
where $f(X) = 8X^3 - 6PX - U = a^{-1}g(2X+b,-3a)$.  Hence a
solution~$x,y\in K$ exists if and only if~$g$ has a linear factor
(with $(2x+b)/(-3a)$ a root of~$g(X,1)$), noting that $x$ cannot be a
double root of~$f$ since~$\D\not=0$, so~$f'(x)\not=0$.
\end{proof}

A special case of
Theorem~\ref{thm:equiv1} below, combined with
Proposition~\ref{prop:linear_factor}, is that all cubics with the same
discriminant which have linear factors over~$K$
are~$\SL(2,K)$-equivalent; this is also implied by the following.

\begin{proposition}\label{prop:standard-root}
  Suppose that $g\in\BC(K;\D)$ has a linear factor in~$K[X,Y]$. Then~$g$
  is~$\SL(2,K)$-equivalent to $Y(X^2 - \frac{1}{4}\D Y^2)$.
\end{proposition}
\begin{proof}
Let the linear factor be~$r X + s Y$ with $r, s\in K$, not both zero.
Let $M=\begin{pmatrix}s&-r\\0&s^{-1}\end{pmatrix}$ if $s\not=0$ and
$M=\begin{pmatrix}s&-r\\r^{-1}&0\end{pmatrix}$ otherwise; then
$\det(M)=1$ and $g^M$ has the linear factor $(r X + s Y)^M = Y$.
Transforming $g$ by~$M$ gives a cubic with coefficients~$(0,b,c,d)$
where~$b\not=0$.  Transforming again by~$\diag(b^{-1},b)$ gives
coefficients~$(0,1,c',d')$, and then
by~$\begin{pmatrix}1&0\\-c'/2&1\end{pmatrix}$ gives
  coefficients~$(0,1,0,d'')$. Comparing discriminants, which do not
  change under unimodular transformations, we see that
  $\D=-4d''$ as required.
\end{proof}

\begin{remark}\label{rmk:Cardano}
Cardano's formula for the roots of the cubic expresses them in terms
of the cube root of a quantity in $L=K(\sqrt{-3\D})$. From the
previous proof we see that, when $P\not=0$, the roots of~$g(X,1)$ are
given by
\[
x = - (b + \sqrt[3]{z} + P/\sqrt[3]{z}) / 3a,
\]
where $z = (U + 3a\sqrt{-3\D})/2$. This is essentially Cardano's
formula.  If $w=\sqrt[3]{z}\in L$ with
$w\overline{w}=\sqrt[3]{\NLK(z)}=P$ then we have
$x=-(b+w+\overline{w})/(3a)\in K$.  All three roots are in~$K$ when
$\sqrt{-3}\notin K$ but~$\sqrt{-3}\in L$, so $L=K(\sqrt{-3})$ and
$\Delta\in(\Ks)^2$, for then $L$ contains a primitive roots of
unity~$\zeta_3$ and $\overline{\zeta_3}=\zeta_3^2$, so that replacing
$w$ by $\zeta_3w$ or $\zeta_3^2w$ in the formula again gives an
element of~$K$.
\end{remark}

\begin{remark}\label{rmk:etale-algebra}
If $-3\D\in(\Ks)^2$, so that $L\cong K\oplus K$, the Cardano
group~$(\Lsmc)_{N=1}$ is isomorphic to~$\Ksmc$.
\end{remark}

\begin{remark}\label{rmk:all_cubics}
  The formulas given here may be used to parametrize all cubic
  extensions of~$K$ as follows. (This construction may also be
  obtained from the Delone-Faddeev parametrization of cubic orders by
  binary cubic forms in~\cite{delone-faddeev}, replacing their base
  ring~$\Z$ by a field.) Each such extension~$M=K(\alpha)$ has a
  discriminant~$\D\in\Ks$, well-defined modulo squares as the
  discriminant of the irreducible cubic minimal polynomial of~$\alpha$
  in~$K[X]$, and $\Delta=1$ (modulo squares) if and only if the
  extension $M/K$ is Galois. To each~$\D\in\Ks$ we form the cubic
  algebra~$L=K[T]/(T^2+3\D)$ as above; then the cubics with
  discriminant~$\D$ (modulo squares) are parametrized by the subgroups
  of the Cardano group~$(\Lsmc)_{N=1}$. Explicitly,
  for~$z\in(\Lsmc)_{N=1}$ with $\NLK(z)=P^3$ and $\TLK(z)=U$ the
  associated cubic is $f_z(X)=X^3-3PX-U$.  (Note that $\overline{z}$
  has the same trace and norm as~$z$, and hence
  $f_{\overline{z}}=f_z$, in accordance with the fact that
  $z\overline{z}=P^3$, so that $\overline{z}$ generates the same
  subgroup of the Cardano group as~$z$.)  The root(s) of $f_z(X)$ are
  $\alpha=w+P/w$ where $w^3=z$.

  A refinement of this construction may be used when $K$ is a number
  field and $S$ a finite set of primes of~$K$, to determine all cubic
  extensions of $K$ which are unramified outside~$S$. For simplicity
  we assume that~$6$ is an~$S$-unit; otherwise the cubic extensions
  constructed may be ramified at primes dividing~$6$ which are not
  in~$S$.  Recall that for an integer $m\ge2$ the subgroup
  \[
  K(S,m) =
  \{x\in\Ks/(\Ks)^m\mid\ord_{\p}(x)\equiv0\pmod{m}\;\forall\p\notin S\}
  \]
  of~$\Ks/(\Ks)^m$ is finite, and may be computed in terms of the
  class group and unit group of~$K$.  First, we determine the possible
  quadratic subfields~$K(\sqrt{\D})$ of the normal closure of the
  cubic, which must also be unramified outside~$S$, by
  restricting~$\D$ to lie in~$K(S,2)$.  For each such~$\Delta$, we
  also have that $K(\sqrt{-3\D})/K$ is unramified outside~$S$; we then
  restrict~$z$ to the subgroup $\ker(\NLK: L(S,3) \to K(S,3))$ in the
  construction of the preceding paragraph.
\end{remark}

\section{Equivalence of binary cubics}\label{sec:equivalence}

In this section we consider necessary and sufficient conditions for
two binary cubic forms over~$K$ with the same discriminant~$\D$ to be
$\GL(2,K)$-equivalent.  We first observe that we need only consider
transformations with determinant~$\pm1$, and reduce the problem to
testing $\SL(2,K)$-equivalence.

\begin{proposition}\label{prop:SL2GL2}
  Let $\D\in\Ks$ and $g_1,g_2\in\BC(K;\D)$.
  \begin{enumerate}
  \item
    If $g_1=g_2^M$ with~$M\in\GL(2,K)$, then~$\det(M) = \pm1$.
  \item
    $g_2$ is~$\GL(2,K)$-equivalent to~$g_1$ if and only if it
    is~$\SL(2,K)$-equivalent to either~$g_1(X,Y)$ or~$g_1(X,-Y)$.
  \end{enumerate}
\end{proposition}
\begin{proof}
If $g_1=g_2^M$ with~$\mu=\det(M)\in\Ks$, then $\disc(g_1) =
\mu^2\disc(g_2)$, so $\mu^2=1$.

For the second part, it suffices to see that~$g_1(X,-Y) = g_1(X,Y)^M$
where $M=\diag(-1,1)$ has determinant~$-1$.
\end{proof}

For $i=1,2$, denote the coefficients of~$g_i\in\BC(K;\D)$
by~$a_i,b_i,c_i,d_i$, their seminvariants by~$P_i,U_i$, and their
Cardano covariants and invariants by~$C_i, z_i$.  Since
$\disc(g_1)=\disc(g_2)=\D$, they have the same covariant
algebra~$L=K[\delta]$ where~$\delta^2=-3\D$.

We will give two criteria for the equivalence of a pair of binary
cubic forms with the same discriminant: one in terms of their Cardano
invariants, and a second one in terms of a cubic bicovariant; the
latter will also directly give a matrix (or matrices)~$M\in\SL(2,K)$
transforming one to the other.  Both criteria are similar to criteria
(established in~\cite{CremonaFisherQuartics}) for the equivalence of a
pair of binary quartic forms with the same invariants.

\subsection{Cubic equivalence in terms of equality of Cardano invariants}
\label{subsec:CardanoEquivalence}

We restate parts~(3) and~(4) of Theorem~\ref{thm:main1}.

\begin{theorem} \label{thm:equiv1}
  Let $g_1,g_2\in\BC(K;\D)$, with common resolvent
  algebra~$L=K[\delta]$.  Then
  \begin{enumerate}
  \item $g_2=g_1^M$ with~$\det(M)=+1$ if and only if
    $z(g_1)=z(g_2)$ in $\Lsmc$;
  \item $g_2=g_1^M$ with~$\det(M)=-1$ if and only if
    $z(g_1)=z(g_2)^{-1}$ in $\Lsmc$.
  \end{enumerate}
  Hence $g_1$ and $g_2$ are $\GL(2,K)$-equivalent if and only if
  $z(g_1)$ and~$z(g_2)$ generate the same subgroup of~$\Lsmc$.
\end{theorem}
\begin{proof}
(1) Suppose that $g_2=g_1^M$ with $M\in\SL(2,K)$.  Then $g_2(X_2,Y_2)
  = g_1(X_1,Y_1)$, where $(X_1\ Y_1) = (X_2\ Y_2)M$, and by the
  covariant property of~$C_1$ we have that
\[
C_2(X_2,Y_2) = C_1^M(X_2,Y_2) = C_1(X_1,Y_1).
\]
It is now clear that the unit values taken by~$C_1$ and~$C_2$ are the
same. Explicitly, we similarly
have~$H_2(X_2,Y_2)=H_1^M(X_2,Y_2)=H_1(X_1,Y_1)$; let $x_2,y_2\in K$ be
such that $H_2(x_2,y_2)\not=0$, and set $(x_1\ y_1)=(x_2\ y_2)M$, so
then also $H_1(x_1,y_1)=H_2(x_2,y_2)\not=0$. (This will be true for
all choices of~$(x_2,y_2)\not=(0,0)$ when $-3\D$ is not a square, and
for all but at most two choices, up to scaling, when $-3\D$ is a
square in~$K$.)  Then $C_1(x_1,y_1) = C_2(x_2,y_2) \in\Ls$, and so
$z(g_1)=z(g_2)$.

For the converse, suppose that~$z(g_1)=z(g_2)$. If~$z(g_1)=z(g_2)=1$,
then by Proposition~\ref{prop:linear_factor}, both cubics have linear
factors, and then Proposition~\ref{prop:standard-root} shows that both
are $\SL(2,K)$-equivalent to $Y(X^2 - \frac{1}{4}\D Y^2)$, and hence
to each other.

Now we may assume that neither cubic has a linear factor; in
particular, $a_1,a_2\not=0$.  Then also $U_1,U_2\not=0$, since
otherwise $g_1(-b_1,3a_1)=a_1U_1=0$ and $g_1(X,Y)$ would be divisible
by~$3a_1X+b_1Y$; similarly for~$g_2$.  By taking suitable
$\SL(2,K)$-transforms of $g_1$ and~$g_2$ if necessary, we can assume
that $P_1,P_2\not=0$ also.  Then we may take as representatives of the
Cardano covariants the elements $z_i=(U_i+3\delta a_i)/2 \in \Ls$ for
$i=1,2$, with~$\NLK(z_i)=P_i^3$.

Since $z_i\overline{z_i} = \NLK(z_i) = P_i^3$, the
condition that $z_2/z_1 \in\Lsc$ is equivalent to~$\overline{z_1}z_2
\in\Lsc$. A computation similar to that in the proof of
Proposition~\ref{prop:linear_factor} shows that $\overline{z_1}z_2
=(x+y\delta)^3$ with $\NLK(x+y\delta)=P_1P_2$ implies that
$x$ is a root of~$f(X)$,
where
\[
f(X) = 16X^3 - 12P_1P_2X - (U_1U_2+27a_1a_2\D).
\]

We now reduce to the case~$a_2U_1-a_1U_2=0$.  If this does not hold,
set $B(X,Y) = U_1g_2(X,Y)-a_1G_2(X,Y) = aX^3+bX^2Y+cXY^2+dY^3$, where
$a=a_2U_1-a_1U_2\not=0$.  One may check that $B(4P_1X+b,-3a) = 4aP_1^3
f(X)$.  Hence, since $z(g_1)=z(g_2)$ implies that $f(X)$ has a root
in~$K$, it follows that~$B$ has a linear factor over~$K$.  After an
$\SL(2,K)$-transformation of~$g_2$ we can take this linear factor
to~$Y$, so that $a_2U_1-a_1U_2=0$.  This transform does not
change~$P_1$, so we still have $P_1\not=0$.

Assuming that $a_2U_1-a_1U_2=0$, let $\lambda = a_2/a_1 = U_2/U_1 \in
\Ks$. The seminvariant syzygy then gives $P_2^3=\lambda^2 P_1^3$.
Since $P_1\not=0$, then also~$P_2\not=0$, and we may set $\mu=\lambda
P_1/P_2$, so that $\lambda=\mu^3$ and $P_2/P_1=\mu^2$.  Let
$g_3(X,Y)=g_1(\mu X,\mu^{-1}Y)$, the transform of~$g_1$
by~$\diag(\mu,\mu^{-1})\in\SL(2,K)$; then the seminvariants
$a_3,P_3,U_3$ of~$g_3$ are the same as those of~$g_2$, and we may
check that~$g_2=g_3^M$
with~$M=\begin{pmatrix}1&0\\(b_2-b_3)/3a_3&1\end{pmatrix}\in\SL(2,K)$. Hence
$g_1$ and~$g_2$ are $\SL_2(K)$-equivalent.

(2) Every~$M\in\GL(2,K)$ with~$\det(M)=-1$ can be written as~$M=DM_1$,
where $M_1\in\SL(2,K)$ and~$D=\diag(-1,1)$.  By part~(1), it suffices
(for both implications) to show that~$z(\tilde{g})=z(g_1)^{-1}$, where
$\tilde{g} = g_1^D$. We have $\tilde{g}(X,Y) = g_1(X,-Y)$, which has
covariants~$\tilde{H}(X,Y)=H_1(X,-Y)$ and~$\tilde{G}(X,Y)=-G_1(X,-Y)$,
so its Cardano covariant~$\tilde{C}$ satisfies
\[
  \tilde{C}(X,-Y) = -\overline{C_1}(X,Y),
\]
from which the syzygy implies
\[
C_1(X,Y) \tilde{C}(X,-Y) = (-H_1(X,Y))^3.
\]
Hence $z(\tilde{g}) = z(g_1)^{-1}$ as required.

The last part is now clear, using Proposition~\ref{prop:SL2GL2}(1).
\end{proof}

One consequence of this result, combined with
Proposition~\ref{prop:linear_factor}, is that two irreducible cubics
cannot be equivalent via matrices of both determinants, ~$+1$
and~$-1$; equivalently, only reducible cubics can have automorphisms
with determinant~$-1$. We will return to automorphisms in
Section~\ref{sec:automorphisms}.

From the preceding proof we may extract the following criterion
for~$\SL(2,K)$-equivalence for cubic forms with no linear factor,
complementing Proposition~\ref{prop:standard-root} in the reducible
case.
\begin{corollary}
Let $g_1,g_2\in\BC(K;\D)$ be irreducible. Then $g_2=g_1^M$
with~$M\in\SL(2,K)$ if and only if the polynomial $f(X) = 16X^3 -
12P_1P_2X - (U_1U_2+27a_1a_2\D)$ has a root in~$K$.
\end{corollary}

To complete the proof of Theorem~\ref{thm:main1} we only need to show
that the Cardano invariant map is surjective.
\begin{proposition}
Every $z \in (\Lsmc)_{N=1}$ arises as $z(g)$ for some~$g\in\BC(K;\D)$.
\end{proposition}
\begin{proof}
Let $z=x+y\delta\in \Ls$ be such that~$\NLK(z)=x^2+3\D y^2=P^3$
with~$P\in\Ks$.  Let~$g\in\BC(K)$ have
coefficients~$(2y/3,0,-P/2y,x/6y^2)$ if $y\not=0$
or~$(0,x/P,0,-\D/4x)$ if~$y=0$. In each case one may check that
$\disc(g)=\D$ and that $U(g)=2x$ and~$P(g)=P$, so that $z(g)=z$ as
required.
\end{proof}

\subsection{Cubic equivalence in terms of factors of a cubic bicovariant}
\label{subsec:BicovariantEquivalence}

In the proof of Theorem~\ref{thm:equiv1}, we saw two conditions
equivalent to the $\SL(2,K)$-equivalence of the cubics~$g_1$
and~$g_2$: that a monic cubic (denoted~$f(X)$ in the proof) had a root
in~$K$, or that a third cubic form~$B(X,Y)$ had a linear factor
over~$K$.  However, we only established these conditions under certain
extra conditions, such as the irreducibility of the~$g_i$.  The cubic
form $B(X,Y)$ used in the proof is the specialisation
at~$(X_1,Y_1,X_2,Y_2) = (1,0,X,Y)$
of~$g_2(X_2,Y_2)G_1(X_1,Y_1)-G_2(X_2,Y_2)g_1(X_1,Y_1)$.  We now show
how linear factors of this ``bicovariant'' directly give
matrices~$M\in\SL(2,K)$ such that~$g_1=g_2^M$ (of which there are at
most three), in all cases.  A twist of this bicovariant similarly
produces matrices (again, at most three) of determinant~$-1$
transforming~$g_1$ into~$g_2$.

We continue to use the earlier notation, where~$g_1,g_2\in\BC(K; \D)$
have covariants~$H_i$ and~$G_i$. To the pair~$(g_1,g_2)$ we
associate a~\emph{bicovariant} in~$K[X_1,Y_1,X_2,Y_2]$:
\[
B_{g_1,g_2}(X_1,Y_1,X_2,Y_2) = g_2(X_2,Y_2)G_1(X_1,Y_1) - G_2(X_2,Y_2)g_1(X_1,Y_1).
\]
This is a bihomogeneous polynomial of bidegree~$(3,3)$; that is, it
is homogeneous of degree~$3$ in each pair of variables separately.  It
is also a bicovariant with respect to $\SL(2,K)\times\SL(2,K)$,
meaning that, for~$M_1,M_2\in\SL(2,K)$, we have
\begin{align*}
B_{g_1^{M_1},g_2^{M_2}}(X_1,Y_1,X_2,Y_2) &=
B_{g_1,g_2}(X_1,Y_1,X_2,Y_2)^{(M_1,M_2)} \\&=
B_{g_1,g_2}(X_1',Y_1',X_2',Y_2')
\end{align*}
where $(X_i'\ Y_i') = (X_i\ Y_i)M_i$ for~$i=1,2$.

\begin{lemma}\label{lem:lem1}
$B_{g_1,g_2}$ is divisible by~$X_1Y_2-X_2Y_1$ if and only if $g_2=\pm
  g_1$.
\end{lemma}
\begin{proof}
Clearly, $B_{g,g}$ is divisible by~$X_1Y_2-X_2Y_1$, and so
is~$B_{g,-g}=-B_{g,g}$.

Conversely, if $B_{g_1,g_2}$ is divisible by~$X_1Y_2-X_2Y_1$ then
$B_{g_1,g_2}(X,Y,X,Y)=0$, so $g_2(X,Y)G_1(X,Y)=g_1(X,Y)G_2(X,Y)$. By
the coprimality of~$g_i$ and~$G_i$ for~$i=1,2$ this implies that
$g_2=\lambda g_1$ and $G_2=\lambda G_1$ for some~$\lambda\in\Ks$; but
the first of these equations implies (using covariance of
the~$G_i$) that $G_2=\lambda^3 G_1$; so $\lambda^2=1$.
\end{proof}

We now consider bilinear factors of~$B_{g_1,g_2}$: factors of
degree~$(1,1)$, of the form
\[
L_M = (X_1\ Y_1)M\begin{pmatrix}0&1\\-1&0\end{pmatrix}
\begin{pmatrix}X_2\\Y_2\end{pmatrix} = -s X_1X_2 + r X_1Y_2 - u Y_1X_2 + t Y_1Y_2,
\]
where $M=\begin{pmatrix}r&s\\t&u\end{pmatrix}$ is a nonzero matrix.
The previous lemma concerned the factor~$L_I$ associated to the
identity matrix~$I$.

\begin{lemma}\label{lem:lem2}
$L_M$ is irreducible if and only if $\det(M)\not=0$.
\end{lemma}
\begin{proof}
$(aX_1+bY_1)(cX_2+dY_2) = L_M$ with
  $M=\begin{pmatrix}a\\b\end{pmatrix} \begin{pmatrix}d&-c \end{pmatrix}$,
  so~$M$ has rank~$1$ if and only if~$L_M$ is reducible.
\end{proof}

\begin{lemma}\label{lem:lem3}
If $L_M$ divides~$B_{g_1,g_2}$ for~$g_1,g_2\in\BC(K;\D)$, then
$\det(M)\in\Ksq$ and~$g_2^M=\pm\det(M)^{1/2}g_1$.
\end{lemma}
\begin{proof}
Suppose that $L_M$ divides~$B_{g_1,g_2}$. If~$\det(M)=0$ then by the
previous lemma, $B_{g_1,g_2}$ is divisible by~$(aX_1+bY_1)$ for
some~$a,b\in K$ not both zero.  Then we have
$B_{g_1,g_2}(b,-a,X_2,Y_2) = 0$ identically, so $g_2(X_2,Y_2)G_1(b,-a)
= G_2(X_2,Y_2)g_1(b,-a)$. Since~$g_1$ and~$G_1$ are coprime,
$g_1(b,-a)$ and~$G_1(b,-a)$ are not both zero, so this equation
contradicts the coprimality of~$g_2$ and~$G_2$.

Let~$\mu=\det(M)$. Working over~$\overline{K}$ we have~$M =
\mu^{1/2}M_1$ with~$M_1\in\SL(2,\overline{K})$.  Then
$L_M=\mu^{1/2}L_{M_1}$, so $B_{g_1,g_2}$ is also divisible
(over~$\overline{K}$) by~$L_{M_1}$.  Hence $B_{g_1^{M_1^{-1}},g_2} =
B_{g_1,g_2}^{M_1^{-1},I}$ is divisible by $L_{M_1}^{M_1^{-1},I} = L_I
= X_1Y_2-X_2Y_1$ (over $\overline{K}$), so $g_1=\pm g_2^{M_1}$ by
Lemma~\ref{lem:lem1}.  Hence $g_2^M=\pm\mu^{1/2}g_1$, which implies
that~$\mu^{1/2}\in\Ks$.
\end{proof}

\begin{proposition}\label{prop:11factors}
  For $g_1,g_2\in\BC(K;\D)$:
\begin{enumerate}
\item
  Every bilinear factor of~$B_{g_1,g_2}(X_1,Y_1,X_2,Y_2)$
  in~$K[X_1,Y_1,X_2,Y_2]$ is (up to scaling) of the form~$L_M$
  with~$M\in\SL(2,K)$ satisfying~$g_1=g_2^M$.
\item
  Every bilinear factor of
  \[
  g_2(X_2,Y_2)G_1(X_1,Y_1) + G_2(X_2,Y_2)g_1(X_1,Y_1).
  \]
  in~$K[X_1,Y_1,X_2,Y_2]$ is (up to scaling) of the form
  $s X_1X_2 + r X_1Y_2 - u Y_1X_2 - t Y_1Y_2$, where
  $M=\begin{pmatrix}r&s\\t&u\end{pmatrix}$ satisfies~$g_1=g_2^M$ and~$\det(M)=-1$.
\item
  There are at most six $M\in\GL(2,K)$ such that~$g_1=g_2^M$, one for
  each bilinear factor in~$K[X_1,Y_1,X_2,Y_2]$ of
  \[
  g_2(X_2,Y_2)^2G_1(X_1,Y_1)^2 - G_2(X_2,Y_2)^2g_1(X_1,Y_1)^2
  \]
  (which is equal to~$g_2(X_2,Y_2)^2H_1(X_1,Y_1)^3 -
  H_2(X_2,Y_2)^3g_1(X_1,Y_1)^2$).
\end{enumerate}
\end{proposition}
\begin{proof}
  (1) By the lemma, we may scale each linear factor~$L_M$ so that
  $\det(M)=1$.

  (2) Replace $g_1(X_1,Y_1)$ by $g_1(X_1,-Y_1)$ and~$G_1(X_1,Y_1)$ by
  $G_1(X_1,-Y_1)$, which is minus the cubic covariant
  of~$g_1(X_1,-Y_1)$, and apply (1).

  (3) Immediate from~(1) and~(2); for the second expression, apply the
  covariant syzygy.
\end{proof}

All that remains to establish Theorem~\ref{thm:main2} is to check that
every transform from~$g_2$ to~$g_1$ arises from bilinear factors in
this way.

\begin{theorem}\label{thm:equiv2}
Let~$g_1,g_2\in\BC(K;\D)$. Then $g_1$ and~$g_2$
are~$\SL(2,K)$-equivalent if and only if~$B_{g_1,g_2}$ has a bilinear
factor in~$K[X_1,Y_1,X_2,Y_2]$, and every bilinear factor of
~$B_{g_1,g_2}$ has the form~$L_M$ with~$M\in\SL(2,K)$, where
$g_1=g_2^M$.
\end{theorem}
\begin{proof}
If $g_1=g_2^M$ with $M\in\SL(2,K)$ then
$B_{g_1,g_2}=B_{g_2,g_2}^{M,I}$, which is divisible by
$(X_1Y_1-X_2Y_1)^{M,I}=L_M$.  The converse is
Corollary~\ref{lem:lem3}.
\end{proof}

\begin{remark}
  It follows from Theorem~\ref{thm:equiv2} and
  Proposition~\ref{prop:11factors} that the number of
  matrices~$M\in\GL(2,K)$ with $g_1=g_2^M$ is at most six, with at
  most three with each determinant~$\pm1$. This number, if nonzero, is
  also the number of automorphisms in~$\GL(2,K)$ of $g_1$ (or of
  ~$g_2$); in the next section, we study automorphisms in more detail,
  and will see that the number of automorphisms in~$\SL(2,K)$ is
  either~$1$ or~$3$, the latter occurring if and only if~$\D\in\Ksq$;
  hence $B_{g_1,g_2}$ splits completely into bilinear factors if and
  only if $g_1,g_2$ are~$\SL(2,K)$-equivalent and $\D\in\Ksq$.
\end{remark}

\section{Automorphisms of cubic forms}\label{sec:automorphisms}

For~$g\in\BC(K)$, the automorphism group of~$g$ is the
subgroup~$\AA(g)$ of matrices~$M\in\GL(2,K)$ such that~$g^M=g$.  Since
for~$M=\lambda I$ we have $g^M=\lambda g$, the only\footnote{Under the
  untwisted action of~$\GL(2,K)$, transforming by~$\zeta I$
  where~$\zeta$ is a cube root of unity is also clearly trivial, and
  the maximum number of possible automorphisms over fields containing
  all cube roots of unity is~$18$.}  scalar matrix in~$\AA(g)$ is the
identity.  By Proposition~\ref{prop:SL2GL2}, every~$M\in\AA(g)$
satisfies~$\det(M)=\pm1$.

Each~$M\in\AA(g)$ acts on the roots of~$g$, viewed as lying
in~$\P^1(\overline{K})$, by linear fractional transformations.  Since
$\PGL(2,\overline{K})$ acts faithfully and transitively on ordered
triples of distinct points on the projective line, the action on the
roots induces an injective homomorphism $\AA(g)\longrightarrow S_3$.
Hence~$\AA(g)$ has order at most~$6$.  Set~$\AA_1(g) =
\AA(g)\cap\SL(2,K)$.  Proposition~\ref{prop:auto2} below implies that
the above homomorphism restricts to an injection of~$\AA_1(g)$ into
the alternating group~$A_3$, so~$\AA_1(g)$ is either trivial or has
order~$3$.

Also, since an automorphism~$M$ has finite order, it has distinct
eigenvalues either both in~$K$ or in a quadratic extension, so $M$ has
exactly two fixed points in~$\P^1(\overline{K})$.

We now determine when~$g$ has automorphisms of order~$2$ and~$3$
respectively.

\begin{proposition}\label{prop:auto2}
$g\in\BC(K)$ has an automorphism~$M$ of order~$2$ if and only if $g$
  has a root in~$\P^1(K)$; then $\det(M)=-1$ and $M$ fixes the root.
\end{proposition}

\begin{proof}
If $M$ is an automorphism of~$g$ of order~$2$, then~$M$ fixes exactly
one root of~$g$.  Moreover, $M^2=I$ but~$M\not=\pm I$, so the
characteristic polynomial of~$M$ is~$X^2-1$; hence $\det(M)=-1$, and
the two eigenvalues~$\pm1$ of~$M$ lie in~$K$.  This implies that the
two fixed points of~$M$ in~$\P^1(\overline{K})$ are $K$-rational; the
root of~$g$ which is fixed by~$M$ is one of these, so is $K$-rational.

Conversely, if $g$ has a $K$-rational root, then without loss of
generality (by Proposition~\ref{prop:standard-root}), $g = Y(X^2 -
\frac{1}{4}\disc(g)Y^2)$, which has the automorphism~$M=\diag(1,-1)$
of order~$2$.
\end{proof}

\begin{proposition}\label{prop:auto3}
  $g\in\BC(K;\D)$ has an automorphism~$M$ of order~$3$ if and only if
  $\D\in\Ksq$; then $\det(M)=+1$ and $M$ acts as a $3$-cycle on the
  roots of~$g$.
\end{proposition}

\begin{proof}
Let $M$ be an automorphism of~$g$ of order~$3$. Then~$M$ acts as a
$3$-cycle on the roots~$\alpha,\beta,\gamma$ of~$g$
in~$\P^1(\overline{K})$, and the characteristic polynomial of~$M$ is
$X^2+X+1$, so~$\det(M)=+1$.  Label the roots so that
$M(\alpha)=\beta$, $M(\beta)=\gamma$, and~$M(\gamma)=\alpha$. Since
the entries of~$M$ are in~$K$ this implies that
$K(\alpha)=K(\beta)=K(\gamma)$.  By the Galois theory of cubics it
follows that $\disc(g)\in\Ksq$; the splitting field of~$g$ is
either~$K$ itself, or a cyclic extension of~$K$.

Conversely, suppose that $\disc(g)\in\Ksq$; there is a unique element
of~$\PGL(2,\overline{K})$ which
maps~$\alpha\mapsto\beta\mapsto\gamma\mapsto\alpha$; by symmetry,
since every polynomial expression in~$\alpha,\beta,\gamma$ which is
fixed by cyclic permutations lies in~$K$, this element actually lies
in~$\PGL(2,K)$.  By a similar argument to that used in
Lemma~\ref{lem:lem3}, any lift to~$\GL(2,K)$ has square determinant,
so there is a lift to~$M\in\SL(2,K)$ such that~$g^M=g$, and $M^3$ acts
trivially on the roots, so~$M^3=I$.
\end{proof}

Putting these parts together, we see that the automorphism group of a
binary cubic form~$g$ depends only on the Galois group of the
splitting field of~$g$.

\begin{theorem}
  Let $g$ be a binary cubic form defined over the field~$K$, with
  $\D\in\Ks$.
  \begin{enumerate}
    \item The group~$\AA_1(g)$ of $\SL(2,K)$-automorphisms of~$g$ is
      trivial unless~$\D$ is a square, in which case it is cyclic of
      order~$3$.
    \item The group~$\AA(g)$ of $\GL(2,K)$-automorphisms of~$g$ modulo
      scalars is isomorphic to a subgroup of the symmetric
      group~$S_3$, namely the centraliser in~$S_3$ of the Galois group
      of~$g$.  Specifically:
      \begin{itemize}
        \item ${\AA}(g)$ is trivial if and only if if~$g$ is
          irreducible over~$K$ and $\D\notin\Ksq$;
        \item ${\AA}(g)\cong C_3$ if and only if~$g$ is
          irreducible over~$K$ and $\D\in\Ksq$;
        \item ${\AA}(g)\cong C_2$ if and only if~$g$ is
          reducible over~$K$ and $\D\notin\Ksq$ (so that~$g$ has
          exactly one root over~$K$);
        \item ${\AA}(g)\cong S_3$ if and only if~$g$ is
          reducible over~$K$ and $\D\in\Ksq$ (so that~$g$ splits
          completely over~$K$).
      \end{itemize}
  \end{enumerate}
\end{theorem}

Part~(2) of this result is the same as Theorem~3.1 in the
work~\cite{xiao2019binary} of~Xiao, which gives the result only
for~$K=\Q$, though the proof Xiao gives is general and similar to
ours. He also states (in Proposition~2.1 of~\cite{xiao2019binary})
that~$\AA(g)\cong S_3$ when $K=\C$.  Automorphisms of binary cubic
forms over~$\Z$ are also the subject of~\cite{neusel1996cubic}, which
also considers the singular case; there, the author's motivation is to
find all subgroups of~$\GL(2,\Z)$ whose invariant subring in~$\Z[X,Y]$
contains a cubic form; they embed~$\AA(g)$ into~$\GL(2,\F_3)$ by
reduction modulo~$3$, and then consider all possible subgroups of that
group, eventually reaching the same conclusion as here.

Part~(2) also follows from the Delone-Faddeev parametrization of cubic
orders by binary cubic forms from~\cite{delone-faddeev}, whose proof
goes through with their base ring~$\Z$ replaced by any field: see
Proposition~12 of~\cite{bhargava-shankar-tsimerman2013}.

\section{Integral equivalence}\label{sec:integral-equivalence}

Our results so far concern binary cubic forms over a field~$K$
(with~$\ch(K)\not=2,3$) and their equivalence under the actions
of~$\SL(2,K)$ and~$\GL(2,K)$.  In applications, one may be interested
in forms with coefficients in some subring~$R$ of~$K$ and their orbits
under~$\SL(2,R)$ and~$\GL(2,R)$.  One classical example is when~$R=\Z$
and~$K=\Q$, or more generally when~$K$ is a number field and~$R$ its
ring of integers.

Using Theorem~18 of~\cite{bhargava-elkies-shnidman2020}, it would be
possible to write down an equivalence test in terms of a generalised
Cardano invariant, at least when~$R$ is a Dedekind domain. However it
is simpler to use our second criterion.  Since $\GL(2,R)$-equivalence
implies~$\GL(2,K)$-equivalence, and similarly for~$\SL(2)$, it is easy
to adapt our results to give a test for~$\GL(2,R)$-{}
or~$\SL(2,R)$-equivalence.  First, let~$g_1,g_2$ be two cubic forms
with coefficients in~$R$ and the same nonzero discriminant~$\D$.
Using Theorem~\ref{thm:equiv2} we can find all~$M\in\SL(2,K)$ such
that~$g_1=g_2^M$, if any.  If there are none, then certainly the forms
are not $\SL(2,R)$-equivalent; if any such~$M$ exist, their number
will be either~$1$ or~$3$, the latter if and only if~$\D$ is a square,
and the forms are $\SL(2,R)$-equivalent if and only if at least one of
the matrices has entries in~$R$.  For $\GL(2,R)$-equivalence, we
repeat using $g_1(X,Y)$ and~$g_2(X,-Y)$.

For example, let~$g_1(X,Y)=X^3-16Y^3$
and~$g_2(X,Y)=g_1(2X,Y/2)=8X^3-2Y^3$.  Both have integral
coefficients, and they are $\SL(2,\Q)$-equivalent
via~$\diag(2,1/2)$. Since the common discriminant is~$-2^83^3$ which
is not a square, this is the unique $\SL(2,\Q)$-equivalence, so they
are not~$\SL(2,\Z)$-equivalent. Since the cubics are irreducible,
neither has an automorphism of determinant~$-1$, hence there are no
other~$\GL(2,\Q)$-equivalences between them and they are therefore
not~$\GL(2,\Z)$-equivalent.

As a second example, let $g(X,Y)=(X-Y)(X-2Y)(X-3Y)$. This has cubic
covariant $G(X,Y)=-2Y(3X-7Y)(3X-5Y)$.  Its self-bicovariant~$B_{g,g}$
has three linear factors:
\begin{align*}
  &Y_1X_2-X_1Y_2,\\
  -3&X_1X_2+7X_1Y_2+5Y_1X_2-13Y_1Y_2,\\
  -3&X_1X_2+5X_1Y_2+7Y_1X_2-13Y_1Y_2.
\end{align*}
From the coefficients of the three factors we find that $\AA_1(g) =
\{I, M, M^2\}$, where $M=\frac{1}{2}\begin{pmatrix}-5&-3\\13&7
\end{pmatrix}$ and $M^2=\frac{1}{2}\begin{pmatrix}7&3\\-13&-5
\end{pmatrix}$. This shows that it is possible for two integral forms to
be~$\SL(2,\Q)$-equivalent via a non-integral matrix~$M$, while also
being~$\SL(2,\Z)$-equivalent. Hence when testing
for~$\SL(2,\Z)$-equivalence (and similarly
for~$\GL(2,\Z)$-equivalence) it is important to consider all possible
$\GL(2,\Q)$-equivalences, to determine whether one is integral.

\section{Applications to elliptic curves}\label{sec:elliptic-curves}
Our work~\cite{CremonaFisherQuartics} with Fisher, on binary quartic
forms and their equivalence, was motivated by the application to
$2$-descent on elliptic curves.  Specifically, there is a bijection
between orbits of binary quartic forms with classical invariants~$I,J$
under a suitably twisted group action, and~$2$-covers of the elliptic
curve~$Y^2=X^3-27IX-27J$, with the $2$-covering maps being given by
the syzygy between covariants of the quartic.

There is a similar application for binary cubic forms, this time to
$3$-isogeny descent on elliptic curves with~$j$-invariant~$0$.  This
connection is treated in detail by Bhargava \textit{et al.}
in~\cite{bhargava-elkies-shnidman2020}.  We briefly describe the
connection here, summarising Section~3
of~\cite{bhargava-elkies-shnidman2020}

Let $E_k$ denote the elliptic curve with Weierstrass
equation~$Y^2=X^3+k$, where $k\in\Ks$.  There is a
$3$-isogeny~$\phi:E_k\to E_{-27k}$ with dual~$\hat{\phi}$.  Elements
of the Galois cohomology group~$H^1(G_K,E_{-27k}[\hat{\phi}])$, which
are locally trivial and so represent~$\phi$-Selmer elements, can be
represented by cubic curves $\CC_g:Z^3=g(X,Y)$,
where~$g\in\BC(K;-108k)$; these are~$\phi$-coverings of~$E_k$.  They
are trivial (coming from $K$-rational points on~$E_k$) if and only
if~$\CC(K)\not=\emptyset$.  The covering map~$\CC_g\to E_k$ is given
by the covariant syzygy~\eqref{eqn:syzygy}, which implies that
\[
(x,y,z) \in \CC_g(K)
\quad\implies \quad
\left(\frac{H(x,y)}{(3z)^2},\frac{G(x,y)}{2(3z)^3}\right)\in E_k(K).
\]
This is stated (with slightly different notation) in Remark~30
in~\cite{bhargava-elkies-shnidman2020}, the syzygy being equation~(16)
there.  In~\cite{bhargava-elkies-shnidman2020} the $\GL(2,K)$-action
is untwisted, but as the results there concern orbits under~$\SL(2,K)$
this is immaterial: for example, ~\cite[Theorem
  27]{bhargava-elkies-shnidman2020} states that there is a bijection
between~$H^1(G_K,E_{-27k}[\hat{\phi}])$ and the set
of~$\SL(2,K)$-orbits on~$\BC(K;-108k)$.

The untwisted~$\GL(2)$-action is also studied by Kulkarni and~Ure
in~\cite{kulkarni-ure2021}, who give (under the assumption that~$K$
contains the cube roots of unity) a different relation between the
(untwisted) $\GL(2,K)$-orbits on $\BC(K)$ and elliptic curves over~$K$
with $j$-invariant~$0$.

\providecommand{\bysame}{\leavevmode\hbox to3em{\hrulefill}\thinspace}
\providecommand{\MR}{\relax\ifhmode\unskip\space\fi MR }
\providecommand{\MRhref}[2]{\href{http://www.ams.org/mathscinet-getitem?mr=#1}{#2}}
\providecommand{\href}[2]{#2}

\end{document}